\documentclass{ifacconf}

\usepackage{natbib}        

\usepackage[dvipdfm]{graphicx} 
\usepackage{subfigure}
\usepackage{amsmath}
\usepackage{amsfonts}
\usepackage[psamsfonts]{amssymb}
\usepackage{color}
\usepackage{comment}
\usepackage{times} 

\newtheorem{lemma}{Lemma}

\newtheorem{assumption}{Assumption}
\newtheorem{remark}{Remark}

\usepackage{epsfig} 

\newcommand{\be}{\begin{equation}}
\newcommand{\ee}{\end{equation}}
\newcommand{\bea}{\begin{eqnarray}}
\newcommand{\eea}{\end{eqnarray}}
\newcommand{\ba}{\begin{array}}
\newcommand{\ea}{\end{array}}
\newcommand{\beas}{\begin{eqnarray*}}
\newcommand{\eeas}{\end{eqnarray*}}
\newcommand{\leftm}{\left[\begin{array}}
\newcommand{\rightm}{\end{array}\right]}
\newcommand{\reals}{\mbox{$\mathbb R$}}

\begin{document}
	\begin{frontmatter}
		
		\title{Congestion management \\ in traffic-light intersections \\ via Infinitesimal Perturbation Analysis\thanksref{footnoteinfo}}
		
		\thanks[footnoteinfo]{Research supported in
part by  NSF under Grant CNS-1239225.}
		
		\author[First]{Carla Seatzu}
		\author[Second]{Yorai Wardi}

		\address[First]{Department of Electrical and Electronic Engineering, University of Cagliari, Italy (e-mail: seatzu@diee.unica.it).}
		\address[Second]{School of Electrical and Computer Engineering, Georgia Institute of Technology, Atlanta, USA (e-mail: ywardi@ece.gatech.edu)}

		\begin{abstract}                
	We present a flow-control technique in traffic-light intersections, aiming at regulating queue lengths
	to given  reference setpoints. The technique is based on multivariable integrators with adaptive gains, computed
	at each control cycle by assessing the IPA gradients of the plant functions. Moreover, the IPA gradients are
	computable on-line despite the absence of detailed models of the traffic flows. The technique is applied to
	a two-intersection system where it exhibits  robustness with respect
	to modeling uncertainties and computing errors, thereby permitting us to   simplify the on-line computations perhaps at the expense
	of accuracy while achieving the desired tracking. We compare, by simulation, the performance of a centralized, joint two-intersection
	control with distributed control of each intersection separately, and show similar performance of the two
	control schemes for a range of parameters.
			\end{abstract}
		
		
	\end{frontmatter}

\section{Introduction}
Infinitesimal Perturbation Analysis (IPA) has been
established as a sample-based technique for sensitivity analysis of Discrete Event Dynamic Systems (DEDS).
Specifically, it gives   formulas or algorithms for the sample derivatives (gradients)
of performance functions with respect
to structural and control variables.  One of its salient features is the simplicity and computational
efficiency of its algorithms in a class of DEDS  which follow
formal rules for  propagation of perturbations  in a network. Furthermore, the algorithms are based on
the monitoring and observation of sample paths associated with the evolving state  of the system, and whenever these
are measurable in real time, the IPA has a potential in control. For extensive presentations of the IPA technique, please
see \cite{Ho91,Glasserman91,Cassandras99}.

The principal application-domain of IPA has been in queueing networks. However, in recent years there has been
a growing interest in fluid queues and their generalization to a class of Stochastic Hybrid Systems (SHS).
There are three reasons for that: (i)  the algorithms for IPA often are simpler in the SHS setting than in their DEDS equivalent models; (ii)
they often require only an observation of the system's sample paths but not detailed or explicit knowledge
of the underlying probability law and hence may be implementable on-line; and (iii) the IPA derivative estimators are unbiased in a larger class of problems in the SHS framework as compared to the setting of DEDS.
Initial developments of IPA in the
SHS framework were presented in
\cite{Cassandras02}, and more-recent general  results as well as surveys can be found   in \cite{Cassandras10, Yao11}.

We point out that representing a discrete event system by an SHS model may introduce errors in the performance
evaluation, and this point was addressed in the aforementioned papers and references therein via simulation-based
parameter optimization. In most of these experiments the underlying system that provided the sample paths
was a DEDS, but the algorithms for the IPA had been derived from the SHS model. Despite this discrepancy the optimization
techniques converged to minimum points (or local minima), thereby suggesting a degree of robustness of them
with respect to errors in the computed IPA gradients. This point plays a role in the later developments in this paper,
as will be seen in the sequel.

While the principal use of IPA has been in optimization, recently we  considered an alternative application
in  performance
regulation. Specifically, we addressed the problem of  performance tracking   in DEDS with time-varying characteristics. The feedback law we chose is
comprised of an integrator with adjustable gain, whose adaptation is based on the IPA derivative of the plant function with respect
to the control variable. One of the key issues   concerns the robustness of the regulation technique with respect
to modeling uncertainties and errors in computing the IPA derivative (\cite{Wardi15}).
An application to throughput regulation in computer processors (\cite{Almoosa12b}) highlights the importance of this issue
since the controller has to run at short cycles, and therefore, if it is robust, it can be designed for simplicity and
speed at the expense of accuracy.  This approach was  justified by simulation results in \cite{Almoosa12b} as well as
an analysis cited therein.

Recently we considered an application of the IPA-based performance regulation technique to congestion management
in
traffic-light intersections, where the objective was to
regulate the queue length in a given direction to a given reference
(\cite{Wardi14}).  This paper extends the results therein in the following two ways: It considers vector tracking by a MIMO system while all previous results on
IPA-based regulation concerned only SISO systems, and it considers the IPA derivative of a queue length
at a given intersection by a control parameter at another intersection.
In particular we compare, by simulation, the performance of a joint centralized controller for a two-intersection
system with a decentralized control where each intersection is controlled by its own parameter.
The results show similar performance for a range of parameters, which indicates the aforementioned
robustness of IPA-based control and justifies the use of the simpler, decentralized scheme.

Applications of IPA to road-congestion management  have been addressed in
\cite{Fu03,Panayiotou05}, and more recently in
\cite{Geng12,Geng13,Geng14,Fleck14}. Of course the traffic control problem has been amply researched
for decades (see, e.g., \cite{Fleck14} for a survey of techniques and results), and the IPA approach
aims at on-line optimization using sample gradients in conjunction with stochastic approximation.
The approach in this paper shares the principle of on-line sample gradients, but deviates from the above-mentioned approach
in that it concerns performance regulation (tracking) and not optimization.
In particular, it considers the control of queue lengths at traffic lights as a mean of
congestion avoidance, and consequently the control laws that we propose   are different.

The rest of the paper is is structured as follows. Section 2 summarizes  our regulation technique in general terms.
Section 3 sets the traffic control problem and analyzes the IPA derivatives, Section 4 provides simulation results, and Section 5 concludes the paper. All of the proofs are relegated to the appendix.

\section{Regulation Algorithm: Integral Control with Adaptive Gain}

Consider the $n$-dimensional discrete-time
control system shown in Figure~1, where $r\in \reals^n$ is the setpoint input vector, $k=1,2,\ldots$, denotes time, $y_{k}\in \reals^n$ is the output vector, $e_{k}\in \reals^n$ is the error signal vector, and $u_{k}\in \reals^n$ is the input to the plant. Suppose first that the plant is a time-varying, memoryless nonlinearity
of the form
\begin{equation}
y_{k}=G_{k}(u_{k}),
\end{equation}
where $G_{k}:\reals^n \rightarrow \reals^n$, $k=1,2,\ldots$, is called the {\it plant function}.
Given a reference input vector $r$, the purpose of the control system is to ensure that
$\lim_{k\rightarrow\infty}y_{k}=r$. To this end we choose the controller to be a linear system defined as
\begin{equation}
u_{k}\ =\ u_{k-1}+A_{k}e_{k-1},
\end{equation}
where $A_k \in \reals^{n\times n}$, and the error signal is defined as
\begin{equation}
e_{k}=r-y_{k}.
\end{equation}

\begin{figure}[!t]
	\centering
	\includegraphics[width=3in]{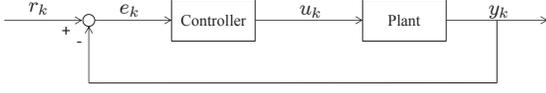}
	\caption{Control System}
	\label{Figure_11}
\end{figure}

Observe that, if $A_k$ is constant independent of $k$, then the above control law essentially is a multi-variable integrator
(adder). Integral controllers often are associated with oscillation and narrow stability margins, therefore we chose a variable-gain integrator to extend the stability margins as well as to guarantee performance of the regulation scheme under
variations in the plant.

We define $A_{k}$ as
\begin{equation}\label{Ak}
A_{k}\ =\left( \frac{\partial G_{k-1}}{\partial u}(u_{k-1})  \right)^{-1} .
\end{equation}

Equations (1) -- (4), computed cyclically in the order
$(4)\rightarrow(2)\rightarrow(1)\rightarrow(3)$, define the dynamics of the closed-loop system.

The rationale behind the choice of $A_k$ in equation~(4) can be seen in the fact that, if the plant is time-invariant and hence $G_k(u) = G(u)$, this control law effectively implements the
Newton-Raphson method for solving the equation $G(u)=r$. This observation was used in \cite{Almoosa12a} for the single-variable control problem of regulating the dynamic power in computer processors. That reference also derived theoretical results concerning robustness of the tracking method to variations in plant modeling and errors in the computation of the gain. Moreover, general results concerning the robustness of the multi-variable Newton-Raphson  method can be found in \cite{Lancaster66}. We will rely on this robustness in the derivation of the IPA gradient in the sequel.

To describe this control system in a temporal framework, when the plant is no more memoryless, let us divide the time axis $\{t\geq 0\}$ into contiguous control cycles, $\Gamma_{1},\Gamma_{2},\ldots$, in increasing order;  $\Gamma_{1}$ starts at time $t=0$, and for every $k=2,3,\ldots$,
$\Gamma_{k}$ starts at the same time $\Gamma_{k-1}$ ends.
Suppose that the quantities $u_{k-1},\ y_{k-1}$, and
$e_{k-1}$ have been computed or derived by the starting time of $\Gamma_{k}$, and
the following sequence of
operations takes place during $\Gamma_{k}$:
\begin{itemize}
	\item[(i)] $A_{k}$ has been computed via Equation (4) during the previous control cycle, and it is available at the starting time of $\Gamma_{k}$;
	\item[(ii)] $u_{k}$ is computed by the controller at the start of $\Lambda_k$ via Equation (2), and the computation is assumed to be instantaneous;  \item[(iii)] the plant acts on the input $u_{k}$ during $\Gamma_{k}$, at the end of which it yields the corresponding output $y_{k}$ via Eq. (1); and \item[(iv)] $e_{k}$ is computed instantaneously by (3) at the end of $\Gamma_{k}$.
\end{itemize}

The most time-consuming  computation can be expected to be  that of $A_{k}$ in Equation (4) since it involves the Jacobian matrix $\frac{\partial G_{k-1}}{\partial u}(u_{k-1}) $,
implicitly assumed to be nonsingular. In the scenario discussed in this paper this is computed by the IPA derivative of the plant function, where we will make approximations designed to simplify the computations while relying 
on the aforementioned robustness. This will be discussed in detail in the next section.
Due to the stochastic, dynamic, and  time-varying  nature of the system it cannot be expected to achieve a perfect performance tracking of a given reference.
Instead, if the regulation algorithm converges faster
than the rate of change of the system, we can expect the performance to chase the desired value in the sense that it approaches it rapidly between drastic changes.
How well this works will be described in Section 4.

\section{Traffic Regulation at Light Intersections}
This section concerns traffic control on a road with two traffic-light intersections, but the analysis appears to be extendable
to a larger number of lights. At each intersection there is a control
parameter associated with the traffic light as defined below, which  can be used to regulate traffic-buildup (queue length) at the intersection.
However, we also consider simultaneous regulation of the two queues. One of the objectives of this study
is to compare the two control strategies.
Whereas the joint control may be more accurate, the decentralized control of each intersection by its parameter
is simpler.
In the derivations of the IPA algorithms
we make judgement calls about simplifying the computations  whenever we feel it to be expedient,
and we test the results by simulation.
The rest of this section defines the problem and  derives the IPA gradients that are used in the regulation scheme.

\subsection{Problem Definition}

Consider the two-intersection road system shown in Figure 2, where each intersection has a traffic light. Assume, for simplicity of argument, that each light cycle consists of red followed by green and there is no orange light.  Let us focus on traffic in the direction of the
horizontal arrows
in the figure. Assume that traffic arrives at the first intersection according to a stochastic process $\{\alpha_{1}(t)\}$,
and let $\delta_1(t)$ denote the instantaneous rate at which it enters the intersection (defined later).
A
fraction $\phi=\phi(t)\in[0,1]$ of the process
$\{\delta_{1}(t)\}$
proceeds to the second
intersection while the rest follows other directions as indicated by the downward arrow.
Also shown  is the part of the cross traffic at the first intersection that is directed to the second intersection,
represented by a stochastic process $\{\tilde{\alpha}_{2}(t)\}$, which interferes and is multiplexed with the horizontal traffic
crossing the first intersection.  The superposition of the two flows comprises the input-flow process
$\{\alpha_{2}(t)\}$ to the second intersection from the left direction, and thus, $\alpha_2(t)=\phi(t)\delta_1(t)+\tilde{\alpha}_2(t)$.
For simplicity's sake we assume no left turns at the intersections or, that left-bound  traffic has its own turn signal and does not interfere
with the traffic in the direction of the arrows that are shown. Thus, a green light in a given direction at an intersection corresponds to a red light
in the perpendicular directions and vice versa.

\begin{figure}[!t]
	\centering
	\includegraphics[width=3in]{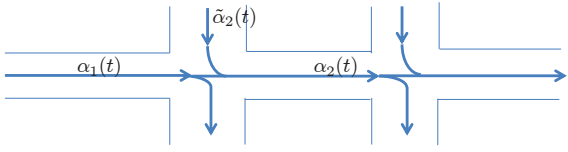}
	\caption{Two-light intersection}
\end{figure}

We consider a scenario in which a High-level (supervisory) controller has computed a traffic plan designed to manage congestion by
balancing  queue lengths,  delays, and light-cycle times at the intersections.
However, traffic bursts and other unpredictable events may cause traffic to deviate from its desirable
behavior  thereby leading to queue buildup. To mitigate these situations we regulate
the intersections' queue lengths by the duty ratios of the light cycles.
We view this procedure as a mechanism for {\it congestion avoidance} which prevents the queue buildup at the second intersection from blocking
traffic at the first intersection. Such acute congestion  is assumed to be handled by a supervisory controller and is
not discussed in this paper.

The regulation technique that we describe is comprised of tracking the queue lengths at the two intersections to given reference values which are assigned,
for example, by
the supervisory controller. The underlying traffic model consists of the fluid-queue system comprised of two queues in
tandem, representing the respective queue buildup at the two intersections in the rightward direction shown in Figure 2.
Let us denote the upstream queue by q1, and the downstream queue by q2. The inflow-rate process to q1 and the outflow
process from it
are $\{\alpha_1(t)\}$ and $\{\delta_1(t)\}$, respectively.
A  $\phi(t)$-fraction of the outflow process  proceeds to the second queue, where it is multiplexed with the
interfering process $\{\tilde{\alpha}_2(t)\}$ to form the inflow-rate process there, $\{\alpha_2(t)\}$,
defined as $\alpha_2(t)=\phi(t)\delta_1(t)+\tilde{\alpha}_2(t)$.

Each of the queues has a constant light cycle of
a given length, $C_1$ and $C_2$, comprised of red followed by green. The control variables
of the system are the lengths of the
red periods at the intersections, denoted by $\theta_{1}$ and $\theta_{2}$, respectively, and we define
$\theta:=(\theta_{1},\theta_{2})^{\top}\in\reals^2$ as the control vector. Naturally the length (duration) of the green periods are $C_i-\theta_i$, $i=1,2$. Note that some of the aforementioned  rate processes depend on $\theta_1$ (like $\beta_1$, $\delta_1$, and $\alpha_2$),
others depend on $\theta_2$ (like $\beta_2$), and $\alpha_1$ and $\tilde{\alpha}_2$ depend on neither $\theta_1$ nor $\theta_2$.
To simplify the notation we will denote their dependence on $\theta$, as in $\beta_1(t,\theta)$, $\alpha_2(t,\theta)$,
etc. The queue-lengths (occupancy) will be denoted by $x_1(t,\theta)$ and $x_2(t,\theta)$, respectively.

During red periods at queue $i$,  $\beta_{i}(t,\theta)=0$.  Upon a light switching from red to green, it  is  realistic
to model the service rate as rising gradually rather than jumping to its highest rate.
The definition of $\beta_{i}(t,\theta)$ reflects this  in the following way: Let $kC_{i}$ be the starting time of the $kth$ cycle at queue $i$, whose red period and green periods
are  the time-intervals
$[kC_{i},kC_{i}+\theta_{i})$ and $[kC_{i}+\theta_{i},kC_{i}+C_{i})$. Let $b_{i}(\tau)$ be a positive-valued, monotone-increasing
random function
of  $\tau\geq 0$. We define $\beta_{i}(t,\theta_{i})=b_{i}(t-(kC_i+\theta_{i}))$ for every $t$ in the green period
$[kC_{i}+\theta_i,kC_{i}+C_{j})$.
Summarizing the definition of $\beta_{i}(t,\theta_{i})$ in both red and green lights, we have,
\begin{eqnarray}
\beta_{i}(t,\theta)=\nonumber \\
\left\{
\begin{array}{ll}
0, & {\rm if}\ t\in [kC_i,kC_{i}+\theta_{i})\\
b_i(t-kC_{i}-\theta_{i}), & {\rm if}\ t\in[kC_{i}+\theta_{i},kC_{i}+C_{i}).

\end{array}
\right.
\end{eqnarray}
Note that this definition of the service rate is quite general, and it includes the special case where
$\beta_{i}(t,\theta)$ holds a constant value $\bar{\beta}_{i}$ during green periods.

Based on the fluid models of the queues, the buffer lengths $x_{i}(t,\theta)$
are related to the inflow  and service rate processes by the following one-sided differential
equation
\begin{equation}
\frac{dx_i}{dt^+}(t,\theta)=\left\{
\begin{array}{ll}
\alpha_{i}(t,\theta)-\beta_{i}(t,\theta), & {\rm if}\ x_{i}(t,\theta)>0\\
0, & {\rm if}\ x_{i}(t,\theta)=0,
\end{array}
\right.
\end{equation}
the outflow rate from the first queue, $\delta_{1}(t,\theta)$, is defined by
\begin{equation}
\delta_{1}(t,\theta)=\left\{
\begin{array}{ll}
\beta_{1}(t,\theta), & {\rm if}\ x_{1}(t,\theta)>0\\\
\alpha_{1}(t), & {\rm otherwise},
\end{array}
\right.
\end{equation}
and the input process to the second queue is defined as
\begin{equation}
\alpha_2(t,\theta)=\phi(t)\delta_{1}(t,\theta)+\tilde{\alpha}_{2}(t).
\end{equation}

The problem that we consider is to regulate the queue lengths by adjusting $\theta$. Specifically, let us divide the
time horizon into contiguous control cycles $\Gamma_k = [\gamma_{k-1}, \gamma_{k})$, $k=1, 2, \ldots$,
and let $T_k:=\gamma_k-\gamma_{k-1}$ be the length of $\Gamma_{k}$.
Defining $y_{i,k}:=\frac{1}{T_k}\int_{\Gamma_k}x_i(t,\theta)dt$
and $y_k:=(y_{1,k},y_{2,k})^{\top}$, the objective is to regulate $y_{k}$ to a given reference vector. The regulation will be carried out via repeated applications of Equations $(4)\rightarrow(2)\rightarrow(1)\rightarrow(3)$ (with $\theta$ replacing $u$)
where the main challenge is
to compute the gain $A_{k}$ defined in (4). This requires the IPA derivative of the plant function $G_{k}(\theta)$ defined as
\begin{equation} \label{G_k}
G_k(\theta)=\displaystyle \frac{1}{T_k}\left[ \begin{array}{c} \displaystyle \int_{\gamma_{k-1}}^{\gamma_k} x_1(t,\theta)dt \\ \displaystyle \int_{\gamma_{k-1}}^{\gamma_k} x_2(t,\theta)dt  \end{array} \right],
\end{equation}
whose computation  is specified in the next subsection.

\subsection{IPA Algorithm}

For the sake of simplicity in the forthcoming discussion we will omit the dependence of $G_{k}(\theta)$, and other terms
in (9), on $k$,   and  denote the generic $k$th control cycle by the interval $\Gamma=[0,T]$.
By Eq.~\eqref{Ak}, it holds that
\begin{equation}
A\ =\left( \frac{\partial G}{\partial \theta}  \right)^{-1}
\end{equation}
where $G(\theta):=\big(G_1(\theta),G_2(\theta)\big)^{\top}$ is given by (9), and we next consider the computation of the IPA derivative
$\frac{\partial G}{\partial \theta}$. Since both $\theta$ and $G(\theta)$ are two-dimensional,
we will be concerned with the four partial derivatives, $\frac{\partial G_{i}}{\partial\theta_{j}}$, for $i,j=1,2.$

\subsubsection{The case where $i=j$}

By Eq.~\eqref{G_k} we have, for $i=1,2$, that
\begin{equation}\label{Gi}
G_i(\theta_i)= \displaystyle \frac{1}{T} \displaystyle \int_{0}^{T} x_i(t,\theta_i)dt,
\end{equation}
and since $x_i(t,\theta_i)$ is continuous in $t$,
\begin{equation}\label{partial_Gi}
\frac{\partial G_i}{\partial \theta_i}= \displaystyle \frac{1}{T} \displaystyle \int_{0}^{T} \frac{\partial x_i}{\partial \theta_i}(t,\theta_i)dt.
\end{equation}
Reference \cite{Wardi14}, considering only a single-queue system, derived the following result for
the term $\frac{\partial x_i}{\partial \theta_i}$.
For $t$ lying in the interior of an empty period in queue  $i$,   $\frac{\partial x_i}{\partial \theta_i}(t,\theta_i)=0$.
On the other hand, for $t$ lying in the interior of a busy period in queue $i$, let $u_{t}$ be the starting time of the busy period containing $t$.
Recall that $kC_{i}$ is the starting time of the $kth$
red period at queue $i$, and let $kC_{i}$, $k=\ell,\ldots,m$, be those points lying in the interval   $[u_{t},t)$.
Then
\begin{eqnarray}\label{partial_xi_CDC}
\frac{\partial x_i}{\partial \theta_i}(t,\theta_i)=\nonumber \\
\sum_{k=\ell}^{m}\beta_i((kC_i)^-, \theta_i)+\beta_i(t, \theta_i)-\beta_i(u_{t}(\theta_i)^+,\theta_i).
\end{eqnarray}
We point out that the rate-terms in \eqref{partial_xi_CDC} can be measured in real time by detecting the speed of passing
vehicles. For details of these derivations please see \cite{Wardi14}.

Eq.~\eqref{partial_xi_CDC} is especially simple when $\beta_i(t,\theta)$ is equal to
a given constant
$\bar{\beta}_{i}>0$ during green-light periods at queue $i$. In that case
it is readily seen that the first sum-term in the RHS of Eq.~\eqref{partial_xi_CDC} is equal to
$(m-\ell+1)\bar{\beta}_{i}$;  $\beta_i(t, \theta_i)\in\{0,\bar{\beta}_{i}\}$ depending on whether
$t$ is in a red period or a green period, respectively; and
similarly $\beta_i(u_{t}(\theta_i)^+,\theta_i)\in\{0,\bar{\beta}_{i}\}$
depending on whether
$t$ is in a red period or a green period, respectively. In this case the computation of
$\frac{\partial G_{i}}{\partial\theta_{i}}$ is a matter of a simple counting process.

\subsubsection{The case where $i\neq j$}

It was mentioned earlier that we do not consider the case where the second queue blocks traffic at the first queue, and therefore
$x_{1}(t,\theta)$ is a function of $\theta_{1}$ but not of $\theta_{2}$. Consequently,
$\frac{\partial G_1}{\partial\theta_{2}}(\theta)=0$.

We next consider the partial derivative term $\frac{\partial x_2}{\partial \theta_1}(t,\theta)$.
We will derive for it an event-based algorithm, where the events in question are light-switchings,
jumps (instantaneous discontinuities) in traffic rates,
and the beginning and end of busy periods at the queues. We say that
two events are independent if neither event causes the other to occur at the same time. The following assumption is quite common
in the literature on IPA of stochastic hybrid systems (e.g., \cite{Cassandras02,Cassandras10}).
\begin{assumption}
	For a given control variable $\theta\in\reals$, w.p.1, no two independent events occur at the same time.
\end{assumption}

We also implicitly assume that all of the derivative terms mentioned in the sequel exist w.p.1, and point out
general and verifyable conditions guaranteeing this assumption (\cite{Cassandras02,Cassandras10}).

Before deriving the term $\frac{\partial x_2}{\partial \theta_1}(t,\theta)$ we mention four types of approximations
that can be practical in implementations. These approximations will be used for computing the IPA derivatives
but not for analysis of the traffic flows themselves which comprise the state of the ``real'' system. Therefore we expect the
regulation scheme to work for a range of errors because of its aforementioned robustness. Some of these approximations are tested via
simulation in Section 4 while others are the subject of current research.
In the first approximation  we assume that instantaneous traffic rates can be measured as, for instance,
by speed detectors that are commonly used in traffic monitoring. Second,  the fraction-term $\phi(t)$ is
a mathematical construct that is well defined only for the fluid-flow model but not for the ``real'', discrete system, and hence we replace it by a term $\phi\in[0,1]$ that can be computed, for example, by averaging traffic flows taken from real-time measurements. 
The derivative term
$\frac{d\phi}{dt}(t)$ is neglected in the computation of
the IPA derivative. Third, as mentioned earlier, we will test the decentralized version of the controller on traffic
obtained from the correlated two-queue system. Fourth, we neglect in the forthcoming analysis the effects of delays
in the control loop, which will be addressed in
the future.

As for the term $\frac{\partial x_2}{\partial \theta_1}(t,\theta)$, consider separately the cases where $t$ lies in the interior of an empty period vs. a busy period in q2.
If $t$ lies in the interior of an empty period in q2 then obviously
\begin{equation}
\frac{\partial x_2}{\partial \theta_1}(t,\theta)=0.
\end{equation}

Consider next the case where  $t$ lies in a busy period in q2, and denote by  $p=p(\theta)$
the starting time of this busy period. Then
\begin{equation}
x_{2}(t)=\int_{p}^t(\alpha_{2}(\tau,\theta)-\beta_{2}(\tau,\theta))d\tau,
\end{equation}
and therefore
\begin{equation}
\frac{\partial x_{2}}{\partial\theta_1}(t,\theta)=
\frac{\partial}{\partial\theta_1}\Big(\int_{p}^t(\alpha_{2}(\tau,\theta)-\beta_{2}(\tau,\theta))d\tau\Big).
\end{equation}
Suppose first that $p$ lies in a q1-green period including the start of such a period, and that $t$ lies in a q1-green
period as well.\footnote{We use the terms ``q1-green period'', ``q1-busy period'', etc. to mean
	a green period in q1, a busy period in q2, etc.}
Let $[\eta_{i},\xi_{i})$, $i=1,\ldots,k$, denote the q1-red periods contained in the interval
$[p,t)$ in increasing order, and define $\xi_{0}:=p$ and $\eta_{k+1}=t$, so that the intervals
$[\xi_{i-1},\eta_{i})$, $i=1,\ldots,k+1$, are q1-green periods contained in the interval $[p,t]$.
\begin{lemma}
	The following equation is in force,
	\begin{eqnarray}
	\frac{\partial x_2}{\partial\theta_1}(t,\theta)\nonumber \\
	=-\big(\alpha_2(p^+,\theta)-\beta_2(p^+,\theta)\big)\frac{\partial p}{\partial\theta_1}\nonumber \\
	+\sum_{i=1}^{k}\big(\alpha_2(\xi_i^-,\theta)-\alpha_2(\xi_i^+,\theta)\big)\nonumber \\
	 +\sum_{i=1}^{k+1}\frac{\partial}{\partial\theta_1}\int_{\xi_{i-1}^+}^{\eta_i^-}\big(\alpha_{2}(\tau,\theta)-\beta_{2}(\tau,\theta)\big)d\tau.
	\end{eqnarray}
	Moreover, the first term in the Right-Hand Side (RHS) of (17) is equal to 0 except in the following
	situation: the start of the q2-busy period at time $p$ is triggered by a jump up in $\beta_{1}(s,\theta)$ at the
	same time, $s=p$. In this case $\frac{\partial p}{\partial\theta_{1}}=1$, and the first term in the RHS of (17) is equal to
	$-\big(\alpha_2(p^+,\theta)-\beta_2(p^+,\theta)\big)$.
\end{lemma}

The proof can be found  in the appendix.\hfill$\Box$
\begin{remark}
	If a point $\tau\in[0,T]$ lies in the interior of a q1-red period then $\alpha_{2}(\tau,\theta)=\tilde{\alpha}_2(\tau,\theta_2)$,
	which is independent of $\theta_1$. Therefore, if $q$ lies in the interior of a q1-red period then
	$\frac{\partial}{\partial\theta_1}\int_{q^+}^{\xi_1^-}\big(\alpha_2(\tau,\theta)-\beta_2(\tau,\theta)\big)d\tau=0$. In
	this   case Equation (17) remains the same except that the sum in the last term of its RHS starts at $i=2$ instead of $i=1$. Likewise,
	if $t$ lies in the interior of a q1-red period then the last term of that sum is $i=k$ and not $i=k+1$.
\end{remark}

The first two terms in the RHS of (17) involve flow rates  which are assumed to be measurable.
It remains to assess the third term in (17), which we next are set to do.

Consider the term
\begin{equation}
\frac{\partial}{\partial\theta_1}
\int_{\xi_{i-1}^+}^{\eta_i^-}
\big(\alpha_2(\tau,\theta)-
\beta_2(\tau,\theta)\big)d\tau,
\end{equation}
and recall that the interval $[\xi_{i-1},\eta_i)$ is a green period in q1.
Let us partition this interval into busy periods and empty periods in q1. Thus, define
$[\tau_{b,1},\tau_{e,1})$ as the first q1-busy period in the interval
$[\xi_{i-1},\eta_{i})$, $[\tau_{e,1},\tau_{b,2})$ is the following empty period,
followed by the next busy period $[\tau_{b,2},\tau_{e,2})$, etc, and let $\tau_{e,m}$ denote the end-point of the
last q1-busy period in the interval
$[\xi_{i-1},\eta_{i}]$. Observe that if $\xi_{i-1}$ lies in a q1-busy period (or is the starting time of such a period) then
$\tau_{b,1}=\xi_{i-1}$, and if $\xi_{i-1}$ lies in a q1-empty period then
$\tau_{b,1}>\xi_{i-1}$. Likewise,
if $\eta_{i}$ lies in a q1-busy period then $\tau_{e,m}=\eta_{i}$, and if $t$ lies in a q1-empty period then
$\tau_{e,m}<\eta_{i}$.
\begin{lemma}
	(I). If the interval $[\xi_{i-1},\eta_{i}]$ is contained in a single q1-busy period then
	\begin{eqnarray}
	\frac{\partial}{\partial\theta_1}
	\int_{\xi_{i-1}^+}^{\eta_{i}^-}
	\big(\alpha_{2}(\tau,\theta)-
	\beta_{2}(\tau,\theta)\big)d\tau \nonumber \\
	=
	\phi\big(\beta_{1}(\xi_{i-1}^{+},\theta)-
	\beta_{1}(\eta_1^-,\theta)\big).
	\end{eqnarray}
	(II). Suppose that the interval $[\xi_{i-1},\eta_{i}]$ is not contained in a single q1-busy period.
	(i). If both  $\xi_{i-1}$ and $\eta_{i}$ are included in two different q1-busy periods then
	\begin{eqnarray}
	\frac{\partial}{\partial\theta_1}\int_{\xi_{i-1}^+}^{\eta_{i}^-}\big(\alpha_2(t,\theta)-\beta_2(t,\theta)\big)dt=\nonumber\\
	\phi\Big(\frac{\partial x_1}{\partial\theta_1}(\xi_{i-1})-
	\big(\alpha_1(\xi_{i-1})-\beta_1
	(\xi_{i-1}^+\big)
	\frac{\partial\xi_{i-1}}
	{\partial\theta_1}\Big)
	\nonumber \\
	+\phi\Big(\beta_1(\tau_{b,m}^+)
	-\beta_1(\eta_1^-)\Big).
	\end{eqnarray}
	(ii). If $\xi_{i-1}$ is contained in the interior of a q1-empty period then the first additive term
	in the RHS os (20) is zero. (iii) If $\eta_{i}$ lies in a q1-empty period then the last additive term
	in the RHS of (20) is zero.
\end{lemma}

The proof can be found  in the appendix.\hfill$\Box$

\begin{remark} In (20), the term $\frac{\partial x_1}{\partial\theta_1}(\xi_{i-1})$  is assumed to be computable via (13) as applied to q1. It holds that
	$\frac{\partial\xi_{i-1}}{\partial\theta_1}=1$  for $i\geq 1$ (since $\xi_{i-1}=kC_1+\theta_1$),  and
	$\frac{\partial\xi_{i-1}}{\partial\theta_1}=\frac{\partial p}{\partial\theta_1}$ for $i=1$ (since $\xi_0=p$).
\end{remark}

The analysis leading to Lemma 1, Lemma 2, and the remarks that follow them imply that Algorithm 1, below,
computes the IPA derivative $\frac{\partial x_{2}}{\partial\theta_1}(t,\theta)$ for $t$ lying in a q2-busy period.
By (9), this will complete the computation of
$\frac{\partial G_2}{\partial\theta_1}(\theta)$. The algorithm's
description
focuses on a single q2-busy period, where we use the notation defined for the analysis in the
earlier paragraphs. Thus, let $p$ be the starting time of the busy period and, in increasing
order,
let  $\eta_{i}$ and
$\xi_{i}$, respectively, $i=1,2,\ldots$, be the starting times of q1-red periods and green
periods during the q1-busy period begun at time $p$. If $p$ lies in a q1-red period then we set $\eta_{1}=p$, while
if $p$ lies in a q1-green period then we define $\xi_{0}=p$ and in this case $\eta_{1}>p$. In either case, $\xi_{1}>p$.

The algorithm computes, recursively, quantities $D_{i}$ at the times $\eta_{i}$, and quantities $E_{i}$ at
times $\xi_{i}$. Furthermore, for every $t\in[\eta_{i},\xi_{i})$ (q1-red period) it sets $\frac{\partial x_{2}}{\partial\theta_1}(t,\theta)=D_i$,
while for every $t\in[\xi_{i},\eta_{i+1})$ (q1-green period) it computes a function
$g_i(t,\theta)$ (defined below) and then sets $\frac{\partial x_{2}}{\partial\theta_1}(t)=E_i+g_i(t)$. The algorithm has the following form.\\

{\bf Algorithm 1}:
\begin{itemize}
	\item
	{\it At time $p$}: \\
	If $p$ lies in  a q1-red period: Set $\eta_{1}=p$,  set $D_{1}=0$. For every $t\in[\eta_{1},\xi_{1})$,
	set $\frac{\partial x_{2}}{\partial\theta_{1}}(t,\theta)=0$.\\
	On the other hand, if $p$ lies in  a q1-green period, set $\xi_{0}=p$, and set $E_0$ as follows:
	If the
	start
	of the busy period is triggered by a
	jump $\beta_{1}(\cdot)$,
	and $p$ lies in a q1-green and busy period  including the start of such a period, set
	\begin{equation}
	E_0=
	-\big(\phi\beta_{1}(p,\theta)-\beta_{2}(p,\theta)\big);
	\end{equation}
	otherwise, set $E_0=0$.
	Moreover, for every $t\in[\xi_0,\eta_1]$, set
	\begin{equation}
	\frac{\partial x_{2}}{\partial\theta_{1}}(t,\theta)=E_{0}+g_{0}(t),
	\end{equation}
	where $g_{0}(t)$ is defined below.
	\item
	{\it At time $\eta_i$, $i=2,\ldots$}:\\
	Set
	\begin{equation}
	D_{i}=E_{i-1}+g_{i-1}(\eta_{i}),
	\end{equation}
	where the function $g_{i-1}(t)$ is specified later.
	Moreover, for every $t\in[\eta_{i},\xi_{i}]$, set
	\begin{equation}
	\frac{\partial x_{2}}{\partial\theta_{1}}(t,\theta)=D_{i}.
	\end{equation}
	Note that the situation  where $i=1$ corresponds to the case where $p$ lies in a q1-red period, which was
	discussed earlier.
	\item
	{\it At time $\xi_{i}$, $i=1,2,\ldots$}:\\
	Set
	\begin{equation}
	E_{i}=D_{i}+\alpha_{2}(\xi_{i}^{-},\theta)-\alpha_{2}(\xi_{i}^{+},\theta).
	\end{equation}
	Moreover, for every $t\in[\xi_{i},\eta_{i+1}]$, define
	\begin{equation}
	\frac{\partial x_{2}}{\partial\theta_{1}}(t,\theta)=E_{i}+g_i(t),
	\end{equation}
	where the function $g_{i}(t)$ is defined as follows.
	\item
	{\it Definition of $g_{i}(t)$ for $t\in [\xi_{i},\eta_{i}]$, $i=0,1,\ldots$:}\\
	If $t$ in the interior of a q1-empty period, set $g_{i}(t)=0$.\\
	On the other hand, if $t$ lies in  a q1-busy period, we set $g_{i}(t)=g_{i,1}(t)+g_{i,2}(t)$, where
	the functions $g_{i,1}$ and $g_{i,2}$ are defined as follows.
	$g_{i,1}(t)$ depends on whether $\xi_{i}$ lies in the interior of a q1-empty period or in a q1-busy period.
	If $\xi_{i}$ lies in the interior of a q1-empty period then $g_{i,1}(t)=0$. On the other hand, if
	$\xi_{i}$ lies in a q1-busy period (including the start of such a period), then we define
	\begin{equation}
	g_{i,1}(t)=\phi\Big[\frac{\partial x_{1}}{\partial\theta_{1}}(\xi_{i},\theta)-\big(\alpha_{1}(\xi_{i},\theta)-\beta_{1}(\xi_{i}^+,\theta)\big)\frac{\partial\xi_{i}}{\partial\theta_{1}}\Big].
	\end{equation}
	We point out that for every $i=1,2,\ldots$, the last term in (28) is
	$\frac{\partial\xi_{i}}{\partial\theta_{1}}=1$. In the special case where $\xi_{0}=p$, we have that
	$\frac{\partial\xi_{0}}{\partial\theta_{1}}=1$ as long as $p$  is a point in a green, busy period at q1, and the start of the q2-busy period begun at time $p$ is due to a jump up in $\beta_{1}$ at time $p$; in all other cases
	$\frac{\partial\xi_{0}}{\partial\theta_{1}}=0$.
	
	Consider next the function $g_{i,2}(t)$. If $t$ lies in the interior of a q1-empty period then $g_{i,2}(t)=0$.
	Likewise, if $\xi_{i}$ lies in a q1-busy period and $t$ lies in the same busy period, then $g_{i,2}(t)=0$.
	For the remaining case $t$ lies in a q1-busy period to which $\xi_{i}$ does not belong. In this case, define $\tau_{i,b}$
	the starting time of the q1-busy period containing $t$. Then, $g_{i,2}(t)$ is defined as
	\begin{equation}
	g_{i,2}(t)=\phi\big(\beta_{1}(\tau_{i,b}^+,\theta)-\beta_{1}(t,\theta)\big).
	\end{equation}
\end{itemize}
\hfill$\Box$
\begin{remark}
	The algorithm looks quite complicated due to the need to track several types of events, but it is quite
	simple to code. Considerable  simplification is likely to result from the likelihood  that a single q1-green
	period is either contained in a single q1-busy period, or divided into a q1-busy period followed by a red period.
	It would be even simpler if the service rates at the queues have the form
	\begin{equation}
	\beta_{i}(t)=\left\{
	\begin{array}{ll}
	0, & {\rm if}\ t\ {\rm lies\ in\ a\ q1-red\ period}\\
	\beta_{i,max}, & {\rm if}\ t\ {\rm lies\ in\ a\ q1-green\ period}.
	\end{array}
	\right.
	\end{equation}
	for given constants $\beta_{i,\max}>0$,
	since in this case $\beta_{i}(t,\theta)$ would be  determined by whether $t$
	lies in a red or green period. Even if this is not the case for ``real'' traffic, such a reasoning can be
	used in the IPA algorithm.
\end{remark}

\section{Simulation Examples}

This section presents simulation examples for testing the effectiveness of the proposed regulation technique.
The traffic-light cycles are is $C_1=C_2=1$, and each  control cycle consists of 20 light cycles.
The process $\{\alpha_1(t)\}$  consists of an off/on model where, in the {\it off} stage $\alpha_1(t)=0$, while for each {\it on}
stage $\alpha_1(t)$ has a single drawn value  uniformly distributed in an interval
$[(1-\zeta)\bar{\alpha}_1,(1+\zeta)\bar{\alpha}_1]$; we chose its mean to be $\bar{\alpha}_1=4.1$, and set
$\zeta=0.3$.
The durations of {\it off} periods and {\it on}
periods are drawn from the uniform distributions
on the intervals $[0,0.02]$ and $[0,0.063]$, respectively.
The process $\{\tilde{\alpha}_2(t)\}$  is generated (drawn) in a similar way except  that its mean value is  $\bar{\alpha}_2=\bar{\alpha}_1/10=0.41$. This is motivated by the fact that the external arrival in the second queue is considered as  noise that is not regulated by the traffic light. The main input flow in the second queue is assumed to come from the first queue, where we took $\phi=0.9$.
The service-rate processes $\{\beta_i(t,\theta)\}$,  satisfy Equation (29) with $\beta_{i,max}=5.0$, $i=1,2$.
The set-point reference vector is $r=[0.1, \; 0.1]^{\top}$, and the initial control variables were set to
$\theta_{1}=\theta_{2}=0.8$.

\begin{figure}[!t]
	\centering
	\includegraphics[width=3.5in]{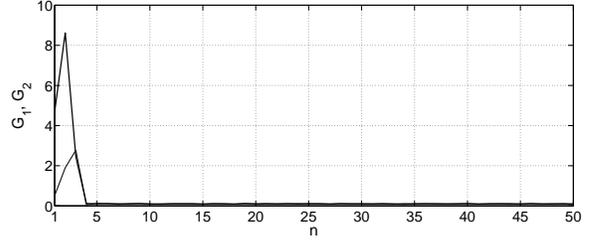}
	\caption{Evolution of $G_{1,k}$ (thick line) and $G_{2,k}$ (thin line) for $k=1\;\ldots,50$}\label{Figure_3}
\end{figure}

\begin{figure}[!t]
	\centering
	\includegraphics[width=3.5in]{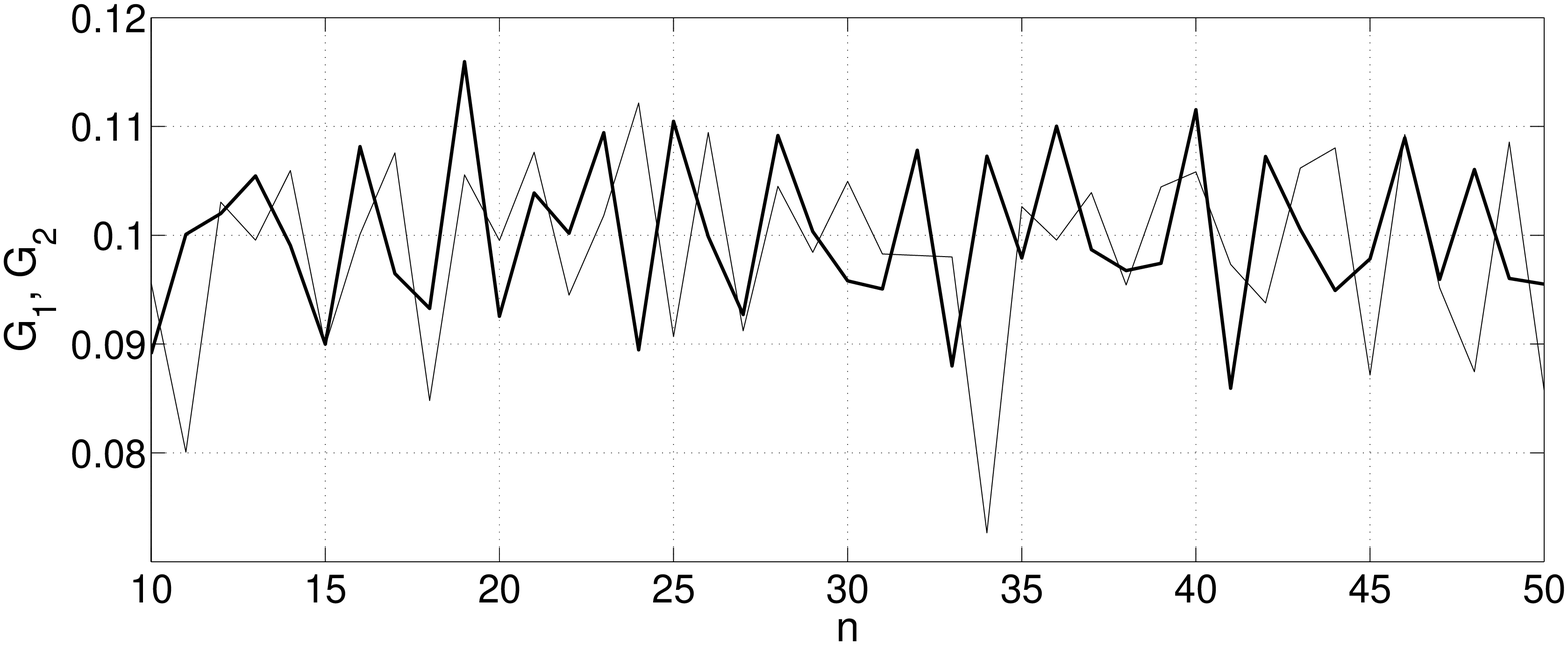}
	\caption{Evolution of $G_{1,k}$ (thick line) and $G_{2,k}$ (thin line) for $k=10\;\ldots,50$}\label{Figure_4}
\end{figure}

Fig.~\ref{Figure_3} depicts the graphs of the obtained outputs $G_{1,k}(\theta)$ and $G_{2,k}(\theta)$
as functions of the counter $k=1,\ldots,50$, and we observe convergence in about 5 iterations.
Fig.~\ref{Figure_4} provides  the same information for $k=10,\ldots,50$ in order to highlight a variability
of the output about the target values of 0.1, which is due to the randomness in the system.  However, the respective means over the last 41 iterations, namely
the quantities $\frac{1}{41}\sum_{k=10}^{50}G_{i,k}(\theta_{k})$, are $0.1001$ and $0.0986$ for $i=1$ and $i=2$, respectively.

\begin{figure}[!t]
	\centering
	\includegraphics[width=3.5in]{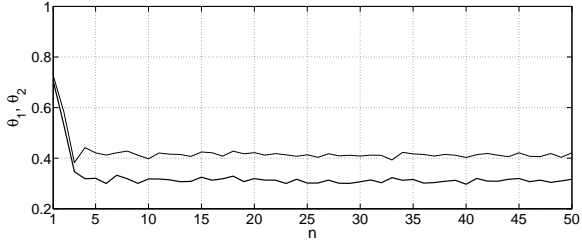}
	\caption{Evolution of $\theta_{1,k}$ (thick line) and $\theta_{2,k}$ (thin line) for $k=1\;\ldots,50$}\label{Figure_5}
\end{figure}

Figure~\ref{Figure_5} shows plots of the control variables $\theta_{i,k}$, $k=1\ldots,50$, and we discern
convergence to their respective values around $\theta_1\sim 0.3113$ and $\theta_2\sim 0.4129$. It is not surprising that the asymptotic value
of $\theta_2$ is larger than that of $\theta_1$. The reason is that the input processes
to the two queues have the same mean rate,  but $\{\alpha_2(t,\theta)\}$ has less variance than
$\{\alpha_1(t,\theta)\}$ due to the action of the first queue. Therefore, to obtain the same mean queue lengths the second queue
would have a larger traffic intensity and hence smaller mean service rate, meaning that $\theta_2>\theta_1$.


A similar behavior, not shown here,  was obtained with different values of the initial control variables as well
as different values of the  traffic parameters  $\zeta$ and $\phi$.

\renewcommand{\baselinestretch}{1.3}

\begin{table*}\label{table1}
	\centering
	\begin{tabular}{|r|c|c|c|c|c|}
		\hline   & $mean|G_1(10:50)-r_1|$ & $mean|G_2(10:50)-r_2|$ & $\max  G_1(1:50)$ & $\max G_2(1:50)$ \\
		
		\hline $\zeta=0.05$, centralized & $0.0843 \cdot 10^{-3}$ & $0.2575\cdot 10^{-3}$ & $0.3507$ & $4.7475$ \\
		decentralized & $0.0855\cdot 10^{-3}$  & $0.0275\cdot 10^{-3}$ &  $0.3481$  & $4.2749$ \\
		
		\hline $\zeta=0.10$, centralized. & $0.0737\cdot 10^{-3}$  & $0.1726\cdot 10^{-3}$ & $0.7258$ & $1.1770$ \\
		decentralized & $0.0637\cdot 10^{-3}$ & $0.0415\cdot 10^{-3}$ & $0.6633$  & $1.9848$ \\
		
		\hline $\zeta=0.15$, centralized. &  $0.1505\cdot 10^{-3}$  & $0.0753\cdot 10^{-3}$ & $2.6368$  & $1.9309$ \\
		decentralized & $0.2167\cdot 10^{-3}$ & $0.0267\cdot 10^{-3}$ &  $2.5817$ & $3.6931$ \\
		
		\hline $\zeta=0.20$, centralized. & $0.1917\cdot 10^{-3}$ & $0.4577\cdot 10^{-3}$  & $4.9831$  & $3.8481$ \\
		decentralized &  $0.2094\cdot 10^{-3}$ & $0.1458\cdot 10^{-3}$ & $5.1314$  & $7.7879$ \\
		
		\hline $\zeta=0.25$, centralized & $0.2642\cdot 10^{-3}$ & $0.6798\cdot 10^{-3}$  & $7.2245$ & $2.2578$ \\
		decentralized & $0.2731\cdot 10^{-3}$ & $7.4621$ & $7.4038$  & $18.3908$ \\
		
		\hline $\zeta=0.30$, centralized & $0.4617\cdot 10^{-3}$  &  $0.2643\cdot 10^{-3}$  & $9.9984$  & $3.0413$ \\
		decentralized & $0.4164\cdot 10^{-3}$ & $8.4869$ &  $10.0046$ & $20.9073$ \\
		\hline
	\end{tabular}
\vspace{0.15in}
	\caption{Comparison between centralized and decentralized regulation scheme.}
\end{table*}

\renewcommand{\baselinestretch}{1}

We compared the joint, two-queue  regulation scheme described in the last two sections with the decentralized
control in which each queue computes its own gain (via (4)) without considering the effects of q1 on q2. The main difference is that
in the second case the Jacobian $\frac{\partial G}{\partial\theta}$ is diagonal, and hence its computation is made
much simpler: compare (13) to Algorithm 1. The comparison was made for various values of the
input variance related to different choices of the traffic parameter $\zeta$. The results, shown in Table 1, comprise averages of
10 independent runs for each indicated value of $\zeta$, where in each run the first central moments of the output
(columns 2 and 3) are computed over
$k=10-50$ in order to avoid the effects of the early transients, while the maximum deviations (columns 4 and 5) are taken over
$k=1-50$.
It is evident that for lower variance $\zeta\leq 0,2$ the performances of the decentralized control is comparable to that of the centralized control. However, as the variance grows the errors associated with neglecting the derivative
$\frac{\partial G_2}{\partial\theta_1}$ render the decentralized control unstable while the centralized controller works well.

\section{conclusions}
The main contribution of this paper is in a flow control
technique in traffic-light intersections which aims at
regulating queue lengths to given reference setpoints.   The technique
is based on multivariable integrators with adaptive gains
computed using IPA. Numerical simulations are presented to
corroborate the effectiveness of the proposed approach. A simpler, decentralized approach based on
approximations designed to reduce the computing efforts is also considered, and its performance is
shown (via simulations) to be similar to that of the multivariable integrators  for a range of parameters.

\section{Appendix}
This section provides proofs to Lemma 1 and Lemma 2.

{\it Proof of Lemma 1.} By (16), Leibnitz rule, and the fact that $\eta_{i+1}=t$ is independent of $\theta_1$,
\begin{eqnarray}
\frac{\partial x_2}{\partial\theta_1}(t,\theta)
=-\big(\alpha_2(p^+,\theta)-\beta_2(p^+,\theta)\big)\frac{\partial p}{\partial\theta_1}\nonumber \\
+\sum_{i=1}^{k}\Big(\alpha_2(\eta_i^-,\theta)-\beta_2(\eta_i^-,\theta)\nonumber \\
-
\big(\alpha_2(\eta_i^+,\theta)-\beta_2(\eta_i^+,\theta_1)\big)\Big)\frac{\partial\eta_i}{\partial\theta_1}\nonumber \\
+\sum_{i=1}^{k}\Big(\alpha_2(\xi_i^-,\theta)-\beta_2(\xi_i^-,\theta)\nonumber \\
-\big(\alpha_2(\xi_i^+,\theta)-\beta_2(\xi_i^+,\theta)\big)\Big)\frac{\partial\xi_i}{\partial\theta_1}\nonumber \\
+\sum_{i=1}^{k+1}\frac{\partial}{\partial\theta_1}\int_{\xi_{i-1}^+}^{\eta_i^-}\big(\alpha_{2}(\tau,\theta)-\beta_{2}(\tau,\theta)\big)d\tau\nonumber \\
+\sum_{i=1}^{k}\frac{\partial}{\partial\theta_1}\int_{\eta_i^+}^{\xi_i^-}\big(\alpha_{2}(\tau,\theta)-\beta_{2}(\tau,\theta)\big)d\tau.
\end{eqnarray}
Light-switchings are events and hence, and by Assumption 1, the functions  $\tilde{\alpha}_2(\tau)$ and
$\beta_2(\tau,\theta)$ are continuous at $\tau=\eta_i$ and $\tau=\xi_i$, $i=1,\ldots,k$.
Next, $\eta_i=kC_1$ and $\xi_i=kC_1+\theta_1$ for some $k=1,\ldots$, hence $\frac{\partial\eta_i}{\partial\theta_1}=0$, and
$\frac{\partial\xi_i}{\partial\theta_1}=1$, rendering the first sum-term in the RHS of (30) to 0, and the last multiplicative
term in the following sum-term to $\frac{\partial\xi_i}{\partial\theta_1}=1$. Furthermore, during q1-red periods $\alpha_2(\tau,\theta)-\beta_2(\tau,\theta)=\tilde{\alpha}_2(\tau)-\beta_2(\tau,\theta_2)$ which is independent of $\theta_1$,
and hence the last additive term in (30) is zero.
In light of all of this, Equation (30) gives (17).

The last assertion of the lemma is obvious.\hfill$\Box$

{\it Proof of Lemma 2.} (I). Suppose that the interval $[\xi_{i-1},\eta_i]$ is contained in a q1-busy period.
By definition, it is also contained  in a q1-green period.
Therefore, for every $\tau$ in this interval,
$\alpha_2(\tau,\theta)=\phi\beta_1(\tau,\theta)+\tilde{\alpha}_2(\tau,\theta)$.
Moreover, by definition and Assumption 1, the functions
$\tilde{\alpha}_2(\tau)$ and $\beta_2(\tau,\theta_2)$ and their jump times are independent of $\theta_1$.
Consequently
\[
\frac{\partial}{\partial\theta_1}
\int_{\xi_{i-1}^+}^{\eta_{i}^-}
\big(\alpha_{2}(\tau,\theta)-
\beta_{2}(\tau,\theta)\big)d\tau=\phi\frac{\partial}{\partial\theta_1}\int_{\xi_{i-1}^-}^{\eta_i^-}\beta_1(\tau,\theta)d\tau.
\]
But for every $\tau\in[\xi_{i-1},\eta_i]$
$\beta_1(\tau,\theta)=b_1(\tau-(kC_1+\theta_1))$ and hence $\frac{\partial\beta_1}{\partial\tau}(\tau,\theta_1)=-\frac{\partial\beta_1}{\partial\theta_1}(\tau,\theta_1)$,
which implies (19).

(II). Consider case (i) where $\xi_{i}$ and $\eta_i$ lie in different q1-busy periods. The times $\tau_{b,j}$ and $\tau_{e,j}$ are event-epoch in q1 and hence, and by Assumption 1, the functions $\beta_2(\tau,\theta)$
and $\tilde{\alpha}_2(\tau)$ are continuous at these points.
Next, $\eta_i$ is a switching time from green to red in q1, hence $\eta_i=kC_{1}$ for some $k=1,\ldots$, and therefore
$\frac{\partial\eta_i}{\partial\theta_1}=0$. Furthermore, recall that the interval
$[\xi_{i-1},\eta_i)$ is contained in a q1-green period. Then for every $\tau\in(\tau_{b,j},\tau_{e,j})$, $j=1,\ldots,m$,
$\alpha_2(\tau,\theta)=\phi\beta_1(\tau,\theta)+\tilde{\alpha}_2(\tau)$,  and for every
$\tau\in(\tau_{e,j},\tau_{b,j+1})$, $j=1,\ldots,m-1$,
$\alpha_2(\tau,\theta)=\phi\alpha_1(\tau)+\tilde{\alpha}_2(\tau)$;
this is independent of $\theta_1$ as is $\beta_2(\tau,\theta)$ and hence
$\frac{\partial}{\partial\theta_1}\big(\alpha_2(\tau,\theta)-\beta_2(\tau,\theta)\big)=0$.
Applying all of this with Leibnitz rule we obtain
\begin{eqnarray}
\frac{\partial}{\partial\theta_1}\int_{\xi_{i-1}^+}^{\eta_i}\big(\alpha_2(\tau,\theta)-\beta_2(\tau,\theta)\big)d\tau\nonumber \\
=\phi\Big[\sum_{j=2}^{m}\big(\alpha_1(\tau_{b,j}^-)-\beta_1(\tau_{b,j}^+,\theta)\big)
\frac{\partial\tau_{b,j}}{\partial\theta_1}\nonumber \\
+\sum_{j=1}^{m-1}\big(\beta_1(\tau_{e,j},\theta)
-\alpha_1(\tau_{e,j})\big)
\frac{\partial\tau_{e,j}}{\partial\theta_1}\nonumber \\
+\sum_{j=1}^m\int_{\tau_{b,j}^+}^{\tau_{e,j}^-}\phi\beta_1(\tau,\theta)d\tau\Big].
\end{eqnarray}
The first sum-term in the RHS of (31) is zero for the following reason:
The point $\tau_{b,j}$ is the starting time of a q1-busy period. If it is triggered by a jump up in $\alpha_1(\tau)$
at $\tau=\tau_{b,j}$ then $\frac{\partial\tau_{b,j}}{\partial\theta_1}=0$ since $\alpha_1(\tau)$ is independent of
$\theta_1$; it cannot be triggered by a jump down in $\beta_1(\tau,\theta)$ since this implies the end of a q1-green period,
but we assume that the interval $[\xi_{i-1},\eta_i)$ is contained in a q1-green period;
and if it is due to a continuous rise of $\alpha_1(\tau)-\beta_1(\tau,\theta)$ from negative to
positive then $\alpha_1(\tau_{b,j})-\beta_1(\tau_{b,j},\theta)=0$.

Next, for every $j=2,\ldots,m-1$, the interval $(\tau_{b,j},\tau_{e,j})$ comprises a q1-busy period, and hence
$\int_{\tau_{b,j}}^{\tau_{e,j}}\big(\alpha_1(\tau)-\beta_1(\tau,\theta)\big)d\tau=0$. Taking derivatives with
respect to $\theta_1$, we obtain,
\begin{equation}
-\frac{\partial}{\partial\theta_1}\int_{\tau_{b,j}^+}^{\tau_{e,j}^-}\beta_1(\tau,\theta)d\tau+\big(\alpha_1(\tau_{e,j})-
\beta_1(\tau_{e,j},\theta)\big)\frac{\partial\tau_{e,j}}{\partial\theta_1}=0.
\end{equation}
The case $j=1$ is different in that
$\int_{\xi_{i-1}}^{\tau_{e,1}}\big(\alpha_1(\tau)-\beta_1(\tau,\theta)\big)d\tau=-x_1(\xi_{i-1},\theta)$, and similarly to
the derivation of (32), we obtain,
\begin{eqnarray}
\big(\alpha_1(\tau_{e,1})-\beta_1(\tau_{e,1},\theta)-\big(\alpha_1(\xi_{i-1}^+)
-\beta_1(\xi_{i-1}^+,\theta)\big)\frac{\partial\xi_{i-1}}{\partial\theta_1}\nonumber \\
=-\frac{\partial x_1}{\partial\theta_1}(\xi_{i-1},\theta).
\end{eqnarray}
By Equation (31) with the aid of (32) and (33), Equation (20) follows.

Finally, parts (II.ii) and (II.iii) of the lemma are obvious in light of the proof of part (II.i), since during a q1-empty period
$x_1(\tau,\theta)=0$ and neither $\alpha_2(\tau,\theta)$ nor $\beta_2(\tau,\theta)$ depend on $\theta_1$.\hfill$\Box$

\bibliography{biblio}

\end{document}